\documentclass[12pt]{article}
\usepackage{amsfonts}
\usepackage{latexsym}
\newtheorem{theorem}{Theorem}
\newtheorem{proposition}{Proposition}
\newtheorem{lemma}{Lemma}
\def\qed{\hfill$\Box$}
\def\grpz{{\mathbb Z}}
\def\CONSTANT{{c_6}}
\begin{document}
\title{A multiplicatively symmetrized version of the Chung-Diaconis-Graham
random process}
\author{Martin Hildebrand
\footnote{Department of Mathematics and Statistics,
University at Albany, State University of New York, Albany, NY 12222.
{\tt mhildebrand@albany.edu}}}
\maketitle
\begin{abstract}
This paper considers random processes of the form $X_{n+1}=a_nX_n+b_n\pmod p$
where $p$ is
odd, $X_0=0$, $(a_0,b_0), (a_1,b_1), (a_2,b_2),...$ are i.i.d., and $a_n$ 
and $b_n$ are independent with $P(a_n=2)=P(a_n=(p+1)/2)=1/2$ and
$P(b_n=1)=P(b_n=0)=P(b_n=-1)=1/3$. This can be viewed
as a multiplicatively symmetrized version of a random process of
Chung, Diaconis, and Graham. This paper shows that order $(\log p)^2$ 
steps suffice
for $X_n$ to be close to uniformly distributed on the integers mod $p$ for
all odd $p$ while order $(\log p)^2$ steps are necessary for $X_n$ to
be close to uniformly distributed on the intgers mod $p$.
\end{abstract}

\section{Introduction}

Chung, Diaconis, and Graham~\cite{cdg} comsidered random processes of
the form $X_{n+1}=2X_n+b_n\pmod p$ where $p$ is odd, $X_0=0$, and
$b_0, b_1, b_2,...$ are i.i.d. with $P(b_n=1)=P(b_n=0)=P(b_n=-1)=1/3$.
They showed that order $(\log p)\log(\log p)$ steps suffice to make
$X_n$ close to uniformly distributed on the integers mod $p$. 
Diaconis~\cite{diaconis} asked about random processes of the form
$X_{n+1}=a_nX_n+b_n \pmod p$ where $p$ is odd, $X_0=0$, and $(a_0,b_0),
(a_1,b_1), (a_2,b_2),...$ are i.i.d. with $a_n$ and $b_n$ being
independent, $P(a_n=2)=P(a_n=(p+1)/2)=1/2$ and $P(b_n=1)=P(b_n=-1)=1/2$.
In his Ph.D. thesis, the author~\cite{mvhphd} showed  that order
$(\log p)^2$ steps suffice to make $X_n$ close to uniformly distributed
on the integers mod $p$ and that order $(\log p)^2$ steps
 are necessary to
make $X_n$ close to uniformly distributed on the integers mod $p$.
The techniques used there can be readily adapted if the distribution 
is changed so that $P(b_n=1)=P(b_n=0)=P(b_n=-1)=1/3$; in this case,
these techniques show that order $((\log p)(\log(\log p)))^2$ steps
suffice to make $X_n$ close to uniformly distributed on the integers mod
$p$ for all odd integers $p$ and order $(\log p)^2$ steps suffice for
almost all odd integers $p$
 while order $(\log p)^2$ steps are necessary to make $X_n$ close
to uniformly distributed in the integrs mod $p$.
This paper shows that this result can be improved to show that
order $(\log p)^2$ steps suffice to make $X_n$ close to uniformly
distributed on the integers mod $p$ for all odd integers $p$.

\section{Some Background, Notation, and Main Result}

We let the integers mod $p$ be denoted by $\grpz/p\grpz$. We may denote
elements of this group by $0, 1,..., p-1$ 
instead of $0+p\grpz, 1+p\grpz,...,(p-1)+\grpz$.

A probability $P$ on the integers mod $p$ satifies $P(s)\ge 0$ for
$s\in\grpz/p\grpz$ and $\sum_{s\in\grpz/p\grpz}P(s)=1$.

We use the variation distance to measure how far a probability
$P$ on $\grpz/p\grpz$ is from the uniform distribution on
$\grpz/p\grpz$. This distance is given by
\[
\|P-U\|=\frac{1}{2}\sum_{s\in\grpz/p\grpz}\left|P(s)-\frac{1}{p}\right|
=\max_{A\subset\grpz/p\grpz}|P(A)-U(A)|
\]
where $P(A)=\sum_{s\in A}P(s)$ and $U(A)=\sum_{s\in A}1/p=|A|/p$.
Note that $\|P-U\|\le 1$ for all probabilities $P$ on $\grpz/p\grpz$.

\begin{proposition}
\label{probmixture}
If $P=p_1P_1+p_2P_2+...+p_mP_m$ where $p_1, p_2,..., p_m$ are positive 
real numbers summing to $1$, then
\[
\|P-U\|\le\sum_{i=1}^mp_i\|P_i-U\|.
\]
\end{proposition}
This proposition can be readily shown using the triangle inequality.

If $P$ is a probability on $\grpz/p\grpz$, define the Fourier tranform
\[\hat P(k)=\sum_{j=0}^{p-1}P(j)e^{2\pi ijk/p}\]
 for $k=0, 1,..., p-1$.
The Upper Bound Lemma of Diaconis and Shahshahani (see, for example,
Diaconis~\cite{diaconis}, p. 24) implies
\[
\|P-U\|^2\le\frac{1}{4}\sum_{k=1}^{p-1}|\hat P(k)|^2.
\]

The main theorem is
\begin{theorem}
\label{mainthm}
Suppose $X_0=0$ and $p$ is an odd integer greater than $1$.
Let $X_{n+1}=a_nX_n+b_n \pmod p$ where $(a_0,b_0), (a_1,b_1), (a_2,b_2),...$
are i.i.d. such that $a_n$ and $b_n$ are independent, 
$P(a_n=2)=P(a_n=(p+1)/2)=1/2$, and $P(b_n=1)=P(b_n=0)=P(b_n=-1)=1/3$.
Let $P_n(j)=P(X_n=j)$ for $j\in\grpz/p\grpz$. Let $\epsilon>0$ be given.
For some $c>0$, if $n>c(\log p)^2$, then $\|P_n-U\|<\epsilon$.
\end{theorem}

\section{Beginnings of the argument}

Observe that 
\begin{eqnarray*}
X_0&=&0\\
X_1&=&b_0\\
X_2&=&a_1b_0+b_1\\
X_3&=&a_2a_1b_0+a_2b_1+b_2\\
&&...\\
X_n&=&a_{n-1}...a_2a_1b_0+a_{n-1}...a_2b_1+...+a_{n-1}b_{n-2}+b_{n-1}
\end{eqnarray*}

We shall focus on the distribution of $X_n$ given values for 
$a_1, a_2,..., a_{n-1}$. In the case where $a_{n-1}=2$, $a_{n-2}=(p+1)/2$,
$a_{n-3}=2$, $a_{n-4}=(p+1)/2$, etc., then
\[
X_n=2(b_{n-2}+b_{n-4}+...)+(b_{n-1}+b_{n-3}+...) \pmod p .
\]
If $n=c(\log p)^2$, then $X_n$ lies between $-(3/2)c(\log p)^2$ and 
$(3/2)c(\log p)^2$ and, for large enough $p$, will not be close to 
uniformly distributed on the integers mod $p$. In the case where
$a_{n-1}=2$, $a_{n-2}=2$, $a_{n-3}=2$, ..., $a_0=2$, then results of
Chung, Diaconis, and Graham~\cite{cdg} show that order $(\log p)\log(\log p)$
steps suffice to make $X_n$ close to uniformly distributed on the integers
mod $p$, and so order $(\log p)^2$ steps suffice as well.

Let $P_n(a_{n-1}, a_{n-2},...,a_1)(s)=P(a_{n-1}...a_1b_0+a_{n-1}...a_2b_1+...+a_{n-1}b_{n-2}+b_{n-1}=s \pmod p )$ where $b_0, b_1,...,b_{n-1}$ are i.i.d.
uniform on $\{1,0,-1\}$.

We shall show
\begin{theorem}
\label{indiviudalcases}
Let $\epsilon>0$ be given.
There exists a constant $c>0$ such that if $n>c(\log p)^2$, then
\[
\|P_n(a_{n-1}, a_{n-2},..., a_1)-U\|<\epsilon/2
\]
except for a set $A$ of values $(a_{n-1}, a_{n-2},..., a_1)$ in 
$\{2,(p+1)/2\}^{n-1}$ where $|A|<(\epsilon/2)2^{n-1}$. ($\{2,(p+1)/2\}^{n-1}$
is the set of $(n-1)$-tuples with entries in $\{2,(p+1)/2\}$.)
\end{theorem}

By Proposition~\ref{probmixture}, Theorem~\ref{indiviudalcases}
implies Theorem~\ref{mainthm}.

\section{Random Walk on the Exponent}
\label{rwexp}

Suppose $a_0, a_1, a_2,...$ are i.i.d. with $P(a_1=2)=P(a_1=(p+1)/2)=1/2$. In
the integers mod $p$, one can view $(p+1)/2$ as $2^{-1}$, the multiplicative
inverse of $2$. So $1, a_{n-1}, a_{n-1}a_{n-2}, a_{n-1}a_{n-2}a_{n-3},...$ can
be viewed as $2^{w_0}, 2^{w_1}, 2^{w_2}, 2^{w_3},...$ where $w_0=0$ and
$w_{j+1}-w_j$ are i.i.d. for $j=0, 1, 2,...$ with 
$P(w_{j+1}-w_j=1)=P(w_{j+1}-w_j=-1)=1/2$. 

Let $M_j=\max\{w_0, w_1,..., w_j\}$ and $m_j=\min\{w_0, w_1,..., w_j\}$.

By Theorem 1 of Section III.7 of Feller~\cite{feller}, 
$P(M_j=\ell)=p_{j,\ell}+p_{j,\ell+1}$ where
$p_{j,\ell}={j \choose (j+\ell)/2}2^{-j}$ where the binomial
coefficient is $0$ unless $(j+\ell)/2$ is an integer between $0$ and $j$,
inclusive. Thus by Central Limit Theorem considerations, for some constant
$c_1>0$, if $\epsilon_1>0$ and $j=\lceil c_1(\log p)^2\rceil$, then
$P(M_j\le 0.5\log_2p)<\epsilon_1/4$ for sufficiently large $p$, and,
by symmetry, 
$P(-m_j\le 0.5\log_2p)<\epsilon_1/4$ for sufficiently large $p$.
Also by Central Limit Theorem considerations, for some constant $c_2>0$,
$P(M_j\ge (c_2/2)\log_2p)<\epsilon_1/4$ and
$P(-m_j\ge (c_2/2)\log_2p)<\epsilon_1/4$ for sufficiently large $p$.
So if
 $j=\lceil c_1(\log p)^2\rceil$, $P(\log_2p<M_j-m_j<c_2\log_2p)>1-\epsilon_1$
for sufficiently large $p$. If this event does not hold, then
$(a_{n-1},a_{n-2},...,a_1)$ might be in the set $A$.

Exercise III.10 of Feller~\cite{feller} gives
\[
z_{r,2n}=\frac{1}{2^{2n-r}}{2n-r\choose n}
\]
where $z_{r,2n}$ is the probability of exactly $r$ returns to the origin
in the first $2n$ steps of the symmetric nearest neighbor random walk on
the integers. Observe
\[
z_{0,2n}=\frac{1}{2^{2n}}{2n\choose n}\sim \frac{1}{\sqrt{\pi n}},
\]
which is approximately a multiple of $1/\log p$ if $n$ is approximately
a multiple of $(\log p)^2$.

Observe that if $r\ge 0$, then
\begin{eqnarray*}
\frac{z_{r+1,2n}}{z_{r,2n}}&=&\frac{1/2^{2n-r-1}}{1/2^{2n-r}}
\frac{{2n-r-1\choose n}}{{2n-r\choose n}}\\
&=&2\frac{n-r}{2n-r}\\
&\le &1.
\end{eqnarray*}
Thus $z_{r+1,2n}\le z_{r,2n}$.

For $k\in[m_j,M_j]$ with $j\lceil c_1(\log p)^2\rceil$,
 let $R(k)$ be the number of $i$
such that $w_i=k$ where $0<i-\min_i\{w_i=k\}\le(\log p)^2$. Observe that
$P(R(k)\le f(p))\le c_3(f(p)+1)/\log p$ for some positive constant
$c_3$.

For some positive constant $c_4$,
observe that 
$E(|\{k:R(k)\le f(p), m_j\le k\le M_j\}|\ | \log_2p<M_j-m_j<c_2(\log_2p))
\le c_4(f(p)+1)$. Thus by Markov's inequality,
$P(|\{k:R(k)\le f(p), m_j\le k\le M_j\}|\ge c_5(f(p)+1)| \log_2p<M_j-m_j<
c_2(\log_2p))\le c_4/c_5$.

\section{Fourier transform argument}

Let $\tilde P_n(a_{n-1}, a_{n-2},..., a_1)(s)=P(2^n(a_{n-1}a_{n-2}...a_1b_0+
a_{n-1}a_{n-2}...a_2b_1+...+a_{n-1}b_{n-2}+b_{n-1})=s \pmod p )$
where $b_0, b_1,..., b_{n-1}$ are i.i.d. uniform on $\{1, 0, -1\}$.
Observe $\|\tilde P_n(a_{n-1}, a_{n-2},..., a_1)-U\|=
\|P_n(a_{n-1}, a_{n-2},..., a_1)-U\|$ since $p$ is odd. Note that all
powers of $2$ in $2^na_{n-1}a_{n-2}...a_1$, $2^na_{n-1}a_{n-2}...a_2$, ...,
$2^na_{n-1}$, $2^n$ are nonnegative.

The Upper Bound Lemma implies
\begin{eqnarray*}
\|\tilde P_n(a_{n-1}, a_{n-2},..., a_1)-U\|&\le&\frac{1}{4}\sum_{m=1}^{p-1}
\prod_{\ell=n+m_j}^{n+M_j}\left(\frac{1}{3}+\frac{2}{3}\cos(2\pi 2^{\ell}m/p)
\right)^{2R(\ell-n)}\\
&\times &\prod_{r=j+1}^{n-1}\left(\frac{1}{3}+\frac{2}{3}\cos(2\pi 2^{n+w_r}m/p)
\right)^2.
\end{eqnarray*}
Note that the first product term is for times up to $j$ and the second product
term is for times after $j$. Recall $j=\lceil c_1(\log p)^2\rceil$.

Note that
\[
\left(\frac{1}{3}+\frac{2}{3}\cos(2\pi 2^{\ell}m/p)\right)^{2R(\ell-n)}
\le
\cases{9^{-R(\ell-n)}&if $1/4\le \{2^{\ell}m/p\}<3/4$\cr
1&otherwise}
\]
and
\[
\left(\frac{1}{3}+\frac{2}{3}\cos(2\pi 2^{n+w_r}m/p)\right)^2
\le
\cases{1/9&if $1/4\le \{2^{n+w_r}m/p\}<3/4$\cr
1&otherwise}
\]
where $\{x\}$ is the fractional part of $x$.

Assume
$|\{k:R(k)\le \CONSTANT \log(\log p), m_j\le k\le M_j\}|<c_5(\log(\log p)+1)$
where $c_5$ is such that $c_4c_6/c_5<\epsilon_2$ where $\epsilon_2>0$
is given and $j=\lceil c_1(\log p)^2\rceil$
and
$|\{k:R(k)<(\log(\log p))^{2.1}, m_j\le k\le M_j\}|<(\log(\log p))^{2.5}$.
Also assume $\log_2p<M_j-m_j<c_2(\log_2p)$,
If these assumptions don't hold, then
 $(a_{n-1},a_{n-2},...,a_1)$ might be in the set $A$.
We shall consider various cases for $m$.

{\underbar {Case 1}}: $m$ is such that for some $\ell\in[n+m_j,n+M_j]$,
$1/4\le\{2^{\ell}m/p\}<3/4$ and $R(\ell-n)>(\log(\log p))^{2.1}$.
Let $S_1$ be the set of such $m$ in $1, 2,..., p-1$. Then, by
arguments similar to those in Chung, Diaconis, and Graham~\cite{cdg}
\[
\sum_{m\in S_1}\prod_{\ell=n+m_j}^{n+M_j}\left(\frac{1}{3}+\frac{2}{3}
\cos(2\pi 2^{\ell}m/p)\right)^{2R(\ell-n)}<\epsilon.
\]
Details appear in Section~\ref{fouriervalues}.

{\underbar {Case 2}}: $m\notin S_1$ and for $b$ values 
of $\ell\in[n+m_j,n+M_j]$, $1/4\le\{2^{\ell}m/p\}<3/4$ and
 $\CONSTANT \log(\log p)<R(\ell-n)\le(\log(\log p))^{2.1}.$ Let $S_{2,b}$ be the
set of such $m$ in $1, 2,..., p-1$.

Let's consider the binary expansion of $m/p$; in particular, consider the
positions $n+m_j+1$ through $n+M_j+1$. If $1/4\le\{2^{\ell}m/p\}<3/4$, then
there is an ``alternation'' between positions $(\ell+1)$ and $(\ell+2)$,
i.e. there is a $1$ followed by a $0$ or a $0$ followed by a $1$. We say
an alternation follows position $\ell$ if there is an alternation between
positions $\ell+1$ and $\ell+2$. Alternations
will start following $b$ of no more than $(\log(\log p))^{2.5}$ positions $\ell$
where $\CONSTANT \log(\log p)<R(\ell-n)<(\log(\log p))^{2.1}$, and alternations
may or may not start following each of no more than
$c_5(\log(\log p)+1)$ positions
$\ell$ with $R(\ell-n)\le \CONSTANT \log(\log p)$.
No other alternations may occur. Place $n+m_j+1$ may be either $0$ or $1$.
 Places $n+m_j+1$ through $n+M_j+1$ of the
binary expansion of $m/p$ are unique for each $m$ in $\{1,2,...,p-1\}$
since $M_j-m_j>\log_2p$ by an observation similar to 
the blocks in the
argument of Chung, Diaconis, and Graham~\cite{cdg} being unique.  So
\begin{eqnarray*}
|S_{2,b}|&\le &2\cdot 2^{c_5(\log(\log p)+1)}{\lfloor(\log(\log p))^{2.5}\rfloor
\choose b}\\
&\le &2\cdot 2^{c_5(\log(\log p)+1)}(\log(\log p))^{2.5b}
\end{eqnarray*}

If $m\in S_{2,b}$, then
\[
\prod_{\ell=n+m_j}^{n+M_j}\left(\frac{1}{3}+\frac{2}{3}\cos(2\pi 2^{\ell}m/p)
\right)^{2R(\ell-n)}
\le
(1/9)^{b\CONSTANT \log(\log p)}.
\]
So
\begin{eqnarray*}
&&\sum_{m\in S_{2,b}}\prod_{\ell=n+m_j}^{n+M_j}\left(\frac{1}{3}+
\frac{2}{3}\cos(2\pi 2^{\ell}m/p)\right)^{2R(\ell-n)}
\\
&\le &
2\cdot 2^{c_5(\log (\log p)+1)}((\log(\log p))^{2.5}
(1/9)^{\CONSTANT \log(\log p)})^b
\end{eqnarray*}

Note that for large enough $p$, $(\log(\log p))^{2.5}(1/9)^{\CONSTANT
\log(\log p)}<1/2$. 
Also observe for $b\ge b_{\min}$ where $b_{\min}$ is a value depending on
$c_5$ and $\CONSTANT$,
\[
2^{c_5(\log (\log p)+1)}((\log(\log p))^{2.5}(1/9)^{\CONSTANT
\log(\log p)})^b\rightarrow 0
\]
as $p\rightarrow\infty$.
Thus
\[
\sum_{b=b_{\min}}^{\infty}
2^{c_5(\log(\log p)+1)}((\log(\log p))^{2.5}
(1/9)^{\CONSTANT \log(\log p)})^b\rightarrow 0
\]
and
\[
\sum_{b=b_{\min}}^{\infty}\sum_{m\in S_{2,b}}\prod_{\ell=n+m_j}^{n+M_j}
\left(\frac{1}{3}+\frac{2}{3}\cos(2\pi 2^{\ell}m/p)\right)^{2R(\ell-n)}
\rightarrow 0.
\]
So all we need to consider are $m\in S_{2,b}$ where $b<b_{\min}$.

To consider such $m$, we shall look at further steps in the Fourier transform.
We shall use the following lemma.
\begin{lemma}
Let $\epsilon^{\prime}>0$ be given. Let $d$ be a positive number. For some
constant $c_7>0$, except with probability no more than $\epsilon^{\prime}$,
\[
\max_{\ell=d+1}^{d+\lfloor c_7(\log p)^2\rfloor}w_{\ell}-\min_{\ell=d+1}^{d+\lfloor c_7(\log p)^2\rfloor}
w_{\ell}>2\log_2p.
\]
If this inequality holds, then, given $m\in\{1, 2,..., p-1\}$, $1/4\le
\{2^{\ell}m/p\}<3/4$ for some $\ell\in\{d+1, d+2,..., d+\lfloor c_7(\log p)^2
\rfloor\}$. With probability at least
$1-(\log(\log p))^{2.5}/\log p$, 
\[
|\{h:\ell+1\le h\le\ell+(\log p)^2, w_{\ell}=w_h\}|>(\log(\log p))^{2.1}.
\]
\end{lemma}

{\it Proof:} Similar to reasoning in section~\ref{rwexp}, the 
existence of $c_7$ follows by Central Limit Theorem 
considerations and Theorem 1 of Section III.7 of Feller~\cite{feller}. 
The existence of such $\ell$ follows since for each positive integer
$k$, at least one of
$\{2^km/p\}$, $\{2^{k+1}m/p\}$,...,$\{2^{k+\lfloor 2\log_2p\rfloor-1}m/p\}$
lies in $[1/4,3/4)$. The result
 on $|\{h:\ell+1\le h\le\ell+(\log p)^2, w_{\ell}=w_h\}|$
follows similarly to the earlier argument 
that $P(R(k)\le f(p))\le c_3(f(p)+1)/\log p$.
\qed

Suppose $n_{before}$ is the number of $m$ being considered, i.e. need further 
Fourier transform terms before going an additional $\lfloor c_7(\log p)^2\rfloor
+\lfloor(\log p)^2\rfloor$ terms. Afterwards, we will need to continue to 
consider only $m$ such that $\ell$ in the lemma exists and 
$|\{h:\ell+1\le h\le \ell+(\log_2p)^2: w_{\ell}=w_jh\}|<(\log(\log p))^{2.1}$;
otherwise we have sufficient additional terms in the Fourier transform; see
Section~\ref{fouriervalues}. 
Except for at most $(\epsilon^{\prime}+o(1))2^{n-1}$ $(n-1)$-tuples in $A$,
$n_{after}\le n_{before}(\log(\log p))^{2.5}/\log p$ where $n_{after}$ is the
number of $m$ still being considered after going the additional
$\lfloor c_7(\log p)^2\rfloor+\lfloor(\log p)^2\rfloor$ steps. 
Repeating this a fixed number $f$
times will give $n_{after}<1$, i.e. $n_{after}=0$ except for at most
$f(\epsilon^{\prime}+o(1))2^{n-1}$ $(n-1)$-tuples in $A$.

\section{Bounding the Fourier transform sums}
\label{fouriervalues}

Some of the ideas in this section, for example ``alternations'', come from
Chung, Diaconis, and Graham~\cite{cdg}.

Suppose $m\in S_1$. If
\[
g(x)=\cases{1/9&if $1/4\le\{x\}<3/4$ \cr 1&otherwise,}
\]
then
\begin{eqnarray*}
\prod_{\ell=n+m_j}^{n+M_j}\left(\frac{1}{3}+\frac{2}{3}\cos(2\pi 2^{\ell}m/p)
\right)^{2R(\ell-n)}
&\le&
\prod_{\ell=n+m_j}^{n+M+j}(g(2^{\ell}m/p))^{R(\ell-n)}\\
&\le&(1/9)^{\CONSTANT \log(\log p)A(B_m)}
\end{eqnarray*}
where $A(B_m)$ is the number of ``alternations'' in the first $M_j-m_j$
positions of the binary expansion of $\{2^{n+m_j}m/p\}$. 
An alternation in the binary expansion
$.\alpha_1\alpha_2\alpha_3...$ occurs when $\alpha_i\ne\alpha_{i+1}$. There
will be an alternation in the first $\lceil\log_2p\rceil$ positions of the
binary expansion of $\{2^{n+m_j}m/p\}$ if $m\in\{1, 2,..., p-1\}$, and for 
different $m\in\{1, 2,..., p-1\}$,  the first $\lceil \log_2p\rceil$
positions of the binary expansion of $\{2^{n+m_j}m/p\}$ will differ.
The inequality ending $<(1/9)^{\CONSTANT \log(\log p)A(B_m)}$ occurs since
for some $\ell\in[n+m_j,n+M_j]$ with $1/4\le\{2^{\ell}m/p\}<3/4$,
$R(\ell-n)\ge(\log(\log p))^{2.1}$ and the $R(\ell-n)$ powers of $1/9$ also
cover all
$c_5(\log(\log p)+1)$ terms of the from $(1/9)^{R(\ell-n)}$ 
with $\ell$ such that
 $R(\ell-n)\le\CONSTANT\log(\log p)$ if $p$ is large enough.

Observe 
\begin{eqnarray*}
\sum_{m\in S_1}(1/9)^{\CONSTANT \log(\log p)A(B_m)}&\le&
\sum_{m=1}^{p-1}(1/9)^{\CONSTANT \log(\log p)A(B_m)}\\
&\le&
2\sum_{s=1}^{M_j-m_j}{M_j-m_j\choose s}(1/9)^{\CONSTANT \log(\log p)s}\\
&\le&
2\sum_{s=1}^{M_j-m_j}(M_j-m_j)^s(1/9)^{\CONSTANT \log(\log p)s}\\
&\rightarrow&0
\end{eqnarray*}
as $p\rightarrow\infty$ if $\log_2p<M_j-m_j<c_2(\log p)$ and $\CONSTANT$
is large enough.

Now suppose $m\in S_{2,0}$ and for some $\ell$ with
$1/4\le\{2^{\ell}m/p\}<3/4$ where $\ell<n-(\log p)^2$ and 
$|\{h:\ell+1\le h\le(\log p)^2,w_{\ell}=w_h\}|\ge(\log(\log p))^{2.1}$, then
\begin{eqnarray*}
&&
\prod_{\ell=n+m_j}^{n+M_j}\left(\frac{1}{3}+\frac{2}{3}\cos(2\pi 2^{\ell}m/p)
\right)^{2R(\ell-n)}
\\
&&
\times
\prod_{r=j+1}^{n-1}\left(\frac{1}{3}+\frac{2}{3}\cos(2\pi 2^{n+w_r}m/p)\right)^2
\\
&\le& (1/9)^{\CONSTANT \log(\log p) A(B_m)}.
\end{eqnarray*}
In other words, the powers of $1/9$ for these 
values of $h$ cover all $c_5(\log(\log p)+1)$
terms of the form $(1/9)^{R(\ell-n)}$ with $\ell$ such that
 $R(\ell-n)\le \CONSTANT\log(\log p)$
if $p$ is large enough. By reasoning similar to the sum involving $m\in S_1$,
\[
\sum_{m\in S_{2,0}}(1/9)^{\CONSTANT \log(\log p) A(B_m)}
\rightarrow 0
\]
as $p\rightarrow\infty$.

\section{Lower Bound}

The argument for the lower bound is more straightforward and is based upon
\cite{mvhphd}.

\begin{theorem}
\label{lowerbound}
Suppose $X_n$, $a_n$, $b_n$, and $p$ are as in Theorem~\ref{mainthm}. Let
$\epsilon>0$ be given. For some $c>0$, if $n<c(\log p)^2$ for large
enough $p$, then $\|P_n-U\|>1-\epsilon$.
\end{theorem}

{\it Proof:}
Let $m_j$ and $M_j$ be as in Section~\ref{rwexp}. For some $c>0$, 
if $n=\lfloor c(\log p)^2\rfloor$, then 
$P(m_j\le -0.25\log_2p)<\epsilon/3$ and $P(M_j\ge 0.25\log_2p)<\epsilon/3$.
If $m_j>-0.25\log_2p$ and $M_j<0.25\log_2p$, then
$2^{\lceil -0.25\log_2p\rceil}X_n$ lies in the interval
$[-{\sqrt p}c(\log p)^2, {\sqrt p}c(\log p)^2]$, and so
$\|P_n-U\|\ge(1-2\epsilon/3)-(2{\sqrt p}c(\log p)^2+1)/p>1-\epsilon$ for
sufficiently large $p$.

\section{Discussion of Generalizations for $a_n$}

One can ask if the results generalize to the case where $a$ is a fixed
integer greater than $1$, $(a,p)=1$, and $P(a_n=a)=P(a_n=a^{-1})=1/2$.
The results indeed should generalize. 
Chapter 3 of Hildebrand~\cite{mvhphd} gives a result if $P(a_n=a)=1$.
This result gives an upper bound similar to the original 
Chung-Diaconis-Graham result with $P(a_n=2)=1$ and involves an $a$-ary
expansion along with a generalization of alternations in a Fourier
transform argument. The random walk on the exponent should work with
powers of $a$ instead of powers of $2$. The Fourier transform
argument may consider the interval $[1/a^2,1-1/a^2)$ instead of
$[1/4,3/4)$. The constant $1/9$ may be replaced by another constant
less than $1$. One needs to be careful with the size of the analogue of
$S_{2,b}$.

Also Breuillard and Varj\'u~\cite{bv} consider the Chung-Diaconis-Graham
process with $P(a_n=a)=1$ where $a$ is not fixed. One might explore
cases where $P(a_n=a)=P(a_n=a^{-1})=1/2$ where $a$ is not fixed but does
have a multiplicative inverse in the integers mod $p$.

\section{Questions for Further Study} 

Eberhard and Varj\'u~\cite{ev} prove and locate
a cut-off phenomonon for most odd integers $p$ in the original
Chung-Diaconis-Graham random process. However, the diffusive
nature of the random walk on the exponent suggests that a cut-off phenomenon
might not appear in the multiplicatively symmetrized version. Exploring
this question more rigorously is a problem for further study.

The Chung-Diaconis-Graham random process can be extended to multiple
dimensions. Klyachko~\cite{klyachko} considers random processes of the
form $X_{N+1}=A_NX_N+B_N \pmod p$ where $X_N$ is a random vector in
$(\grpz/p\grpz)\times(\grpz/p\grpz)$ and $A_N$ is a fixed $2\times 2$
matrix with some conditions. Perhaps techniques 
in this paper could be combined with
Klyachko's result to get a result for the
case where $A_N$ is a fixed $2\times 2$ matrix
or its inverse with probability $1/2$ each.

\section{Acknowledgment}

The author would like to thank the referee for some suggestions. 

This is 
a preprint of an article published in {\it Journal of Theoretical Probability}.
The final authenticated version is available online at

{\tt https://doi.org/10.1007/s10959-021-01088-3}.

\end{document}